\documentclass{amsart}

\usepackage{amssymb}
\usepackage{latexsym}
\usepackage{euscript}
\usepackage{graphics}
\usepackage{amsmath,amsthm}
\usepackage{graphics}
\usepackage[english]{babel}

   \def\sH{{\mathfrak H}}

      \def\dC{{\mathbb C}}

   \def\dN{{\mathbb N}}   
      \def\dR{{\mathbb R}}

   \def\dZ{{\mathbb Z}}

   \def\cW{{\mathcal W}}

\def\e{\varepsilon}
\def\wB{\widetilde{B}}

\def\wb{\widetilde{b}}

\def\f{\varphi}

\def\wt#1{{{\widetilde #1} }}

\def\sg{\operatorname{sign}}

\newcommand{\sddots}{\begin{picture}(2,2)
\multiput(0,0)(1.5,1){3}{.}\end{picture}}

\newtheorem{theorem}{Theorem}[section]
\newtheorem{proposition}[theorem]{Proposition}

\newtheorem{definition}[theorem]{Definition}
\theoremstyle{remark}

\newtheorem{remark}[theorem]{Remark}

\numberwithin{equation}{section}

\title[Inverse  problems]
{Inverse problems for periodic\\ generalized Jacobi matrices}

\author{Maxim Derevyagin}
\address{
Maxim Derevyagin\\
Department of Mathematics MA 4-5\\
Technische Universit\"at Berlin\\
Strasse des 17. Juni 136\\
D-10623 Berlin\\
Germany}
\email{derevyagin.m@gmail.com}

\begin{document}


\begin{abstract}
Some inverse problems for semi-infinite periodic generalized Jacobi matrices are considered.
In particular, a generalization of the Abel criterion is presented. The approach is based on
the fact that the solvability of the Pell-Abel equation is equivalent
to the existence of a certainly normalized $J$-unitary $2\times 2$-matrix polynomial (the monodromy matrix).
\end{abstract}

\maketitle

\section{Introduction}

In 1826, Abel proved that the square root $\sqrt{R}$ of the polynomial $R$
of degree $2n$ without multiple zeros can be expanded into a periodic continued 
fraction, that is  
\begin{equation}\label{defCF}
\sqrt{R}-U=
\frac{1}{{\mathfrak p}_0+\displaystyle{\frac{1}{{\mathfrak p}_1+\displaystyle{\frac{1}{\ddots}}}}}=
\frac{1}{{\mathfrak p}_0}
\begin{array}{l} \\ +\end{array}
\frac{1}{{\mathfrak p}_1}
\begin{array}{l} \\ +\end{array}
\begin{array}{l} \\ \cdots \end{array}
\begin{array}{l} \\ +\end{array}
\frac{1}{{\mathfrak p}_{s-1}}
\begin{array}{l} \\ +\end{array}
\frac{1}{{\mathfrak p}_0}
\begin{array}{l} \\ +\cdots \end{array},
\end{equation}
where $U$ and ${\mathfrak p}_0$, \dots, ${\mathfrak p}_{s-1}$ are polynomials,
 if and only if the Pell-Abel equation
\begin{equation}\label{PellAbeleq}
X^2-RY^2=1
\end{equation}
has the polynomial solutions $X$ and $Y$ ({\it the Abel criterion}). 
Actually, the basic idea of Abel was to find out if there exists a polynomial $\rho$
such that  the integral
\[
\int\frac{\rho(t)}{\sqrt{R(t)}}dt
\]
can be expressed in terms of the elementary functions.
It turned out that such a polynomial $\rho$ exists if and only if the Pell-Abel 
equation~\eqref{PellAbeleq} is solvable in polynomials. The proof of these results
and much more modern information on this topic from the algebraic point of view can be found
in~\cite{Mal01}.

Another appearance of the Pell-Abel equation is intimately related to extremal polynomials.
For example, the Chebyshev polynomials $T_j$ and $U_j$ of the first and second kind satisfy
the following relation
\[
T_j^2(\lambda)-(\lambda^2-1)U_j^2(\lambda)=1.
\]
This observation can be generalized to the case of the polynomials that
least deviate from zero on several intervals~\cite{SYu}. Furthermore, by setting 
\[
R(\lambda)=(\lambda^2-1)\prod_{j=1}^{n-1}(\lambda-\alpha_j)(\lambda-\beta_j),
\]
where $-1<\alpha_1<\beta_1<\dots<\alpha_{n-1}<\beta_{n-1}<1$, we have that 
the equation~\eqref{PellAbeleq} is solvable in polynomials if and only if the set
\[
[-1,1]\setminus\bigcup_{j=1}^{n-1}(\alpha_j,\beta_j)
\]
coincides with the spectrum of a bi-infinite periodic classical Jacobi matrix
~\cite{SYu}. In turn, such a spectrum gives rise to periodic solutions of the Toda
lattice (for example, see~\cite[Chapter~12]{Teschl}).

On the other hand, we have seen semi-infinite classical Jacobi matrices at the beginning
although it was implicitly. Indeed, in some special cases the continued fraction~\eqref{defCF}
can be a $J$-fraction generating a classical Jacobi matrix~\cite{Ach61}. However,
$J$-fractions do not cover all the possible cases of periodic fractions of the form~\eqref{defCF}. 
That is why it is more natural to consider semi-infinite generalized Jacobi matrices associated with continued fractions
of the form~\eqref{defCF}.
Such matrices were introduced in~\cite{DD} and the theory was further developed in~\cite{De09,DD07}.
In particular, direct spectral problems for the periodic generalized Jacobi matrices were considered
in~\cite{De09}. The goal of the present paper is to solve some inverse problems for semi-infinite
generalized Jacobi matrices associated with periodic continued fractions.
It is done in the following way. First, a one-to-one correspondence between certainly normalized 
$2\times 2$-matrix polynomials ({\it the monodromy matrices}) and 
the periodic generalized Jacobi matrices in question is established.
Then we show that the solvability of the Pell-Abel equation is equivalent
to the existence of a normalized $J$-unitary $2\times 2$-matrix polynomial,
that is the monodromy matrix of the underlying periodic generalized Jacobi matrix.
The latter statement generalizes the Abel criterion.
Namely, we give necessary and sufficient conditions for functions of the form
\[
\frac{\sqrt{R(\lambda)}-U(\lambda)}{V(\lambda)}
\]
to be the $m$-functions of semi-infinite periodic generalized Jacobi matrices 
or, equivalently, to admit the periodic continued fraction expansions.

Finally, it should be mentioned that some inverse spectral problems for finite 
generalized Jacobi matrices were studied in~\cite{De06}.

\section{Preliminaries}

\subsection{P-fractions}
Let $\varphi$ be a nonrational function holomorphic at a neighborhood of infinity
and having the property
\[
\overline{\varphi(\overline{\lambda})}=\varphi(\lambda).
\]
So, $\f$ has the following representation in a neighborhood of infinity
\begin{equation}\label{phi_inf}
\displaystyle{\f(\lambda)=-\sum_{j=0}^{\infty}\frac{s_j}{\lambda^{j+1}}},
\end{equation}
where the {\it moments} $s_j$ are real.
A number $n_j\in\dN$ is called
{\it a normal index} 
of the sequence ${\bf s}:=\{s_j\}_{j=0}^\infty$
if $\det (s_{i+k})_{i,k=0}^{n_j-1}\ne 0$.
Since $\varphi$ is not rational, there exists an infinite number
of normal indices of ${\bf s}$ (see~\cite[Section~16.10.2]{G}).
Also, without loss of generality we will always assume that the sequence ${\bf s}$ 
is {\it normalized}, i.e. the first
nonvanishing moment $s_{n_1-1}$ has modulus 1.

As is known~\cite{Mag1}, the series~\eqref{phi_inf} leads to the following infinite continued fraction
\begin{equation}\label{ContF2}
-\frac{\varepsilon_0}{p_0(\lambda)}
\begin{array}{l} \\ - \end{array}
\frac{\varepsilon_0\varepsilon_1b_0^2}{p_1(\lambda)}
\begin{array}{ccc} \\ - & \cdots & -\end{array}
\frac{\varepsilon_{j-1}\varepsilon_j b_{j-1}^2}{p_j(\lambda)}
\begin{array}{cc} \\ - & \cdots \end{array},
\end{equation}
where $\e_j=\pm 1$, $b_j>0$ and
$p_{j}(\lambda)={\lambda}^{k_{j}}+p_{k_{j}-1}^{(j)}{\lambda}^{k_{j}-1}+\dots+
p_{1}^{(j)}\lambda+p_{0}^{(j)}$ are real monic polynomials of
degree ${k_j}$ (see also~\cite{ADL07},~\cite{De}). Note, that
$n_j=k_0+k_1+\dots+k_{j-1}$.

The continued fraction~\eqref{ContF2} is called a $P$-fraction. Actually, the $P$-fraction can be 
considered as a sequence of the linear-fractional transformations~\cite[Section 5.2]{JTh}

\[
T_j(\omega):=\frac{-\e_j}{p_j(\lambda)+\e_jb_j^2\omega}
\]
having the following matrix representation
\begin{equation}\label{Wj}
\cW_j(\lambda)=\begin{pmatrix}0 & -\frac{\varepsilon_j}{b_j}\\
                            \varepsilon_jb_j &  \frac{p_j(\lambda)}{b_j}
                            \end{pmatrix},\quad j\in\dZ_+.
\end{equation}
The superposition $T_{0}\circ T_{1}\circ\dots \circ T_{j}$ of the
linear-fractional transformations corresponds to the product of
the matrices $\cW_l(\lambda)$
\begin{equation}\label{W}
\cW_{[0,j]}(\lambda)=(w_{ik}^{(j)}(\lambda))_{i,k=1}^2:=\cW_0(\lambda)\cW_1(\lambda)\dots
\cW_j(\lambda).
\end{equation}

It is well known that the entries of $\cW_{[0,j]}$ can be expressed in terms of
denominators and numerators of the convergents to the $P$-fraction~\eqref{ContF2}. 
To give these formulas explicitly define the polynomials
$P_{j+1}(\lambda)$, $Q_{j+1}(\lambda)$ by the equalities
\[%\begin{equation}\label{PQj}
   \left(%
\begin{array}{c}
  -Q_0 \\
  P_0 \\
\end{array}%
\right)=\left(%
\begin{array}{c}
  0 \\
  1 \\
\end{array}%
\right),\quad
\left(%
\begin{array}{c}
 -Q_{j+1}(\lambda) \\
 P_{j+1}(\lambda) \\
\end{array}%
\right):=\cW_{[0,j]}(\lambda)\left(%
\begin{array}{c}
  0 \\
  1 \\
\end{array}%
\right), \quad j\in\dZ_+ .
\]%\end{equation}
The relation
$\cW_{[0,j]}(\lambda)=\cW_{[0,j-1]}(\lambda)\cW_{j}(\lambda)$ yields
\[%\begin{equation}\label{W0j}
\cW_{[0,j]}(\lambda)\left(%
\begin{array}{c}
  1 \\
  0 \\
\end{array}%
\right)=
\cW_{[0,j-1]}(\lambda)\left(%
\begin{array}{c}
  0 \\
  \varepsilon_jb_j \\
\end{array}%
\right)=
\left(%
\begin{array}{c}
  -\varepsilon_jb_j Q_j(\lambda) \\
  \varepsilon_jb_j P_j(\lambda) \\
\end{array}%
\right), \quad j\in\dN .
\]%\end{equation}
So, the matrix $\cW_{[0,j]}(\lambda)$ has the form
\begin{equation}\label{SolMatr}
\cW_{[0,j]}(\lambda)=\left(%
\begin{array}{cc}
  -\varepsilon_jb_j Q_j(\lambda) & -Q_{j+1}(\lambda) \\
  \varepsilon_jb_j P_j(\lambda) & P_{j+1}(\lambda) \\
\end{array}%
\right), \quad j\in\dZ_+ .
\end{equation}
Further, the equality
\[%begin{equation}\label{W0j2}
\left(%
\begin{array}{c}
 -Q_{j+1}(\lambda) \\
 P_{j+1}(\lambda) \\
\end{array}%
\right)=\cW_{[0,j-1]}(\lambda)\cW_{j}(\lambda)\left(%
\begin{array}{c}
  0 \\
  1 \\
\end{array}%
\right)=
\frac{1}{b_j}\cW_{[0,j-1]}(\lambda)\left(%
\begin{array}{c}
  -\varepsilon_j \\
  p_j(\lambda) \\
\end{array}%
\right), \quad j\in\dN,
\]%end{equation}
 shows that the polynomials $P_{j}(\lambda)$, $Q_{j}(\lambda)$ are solutions of the difference equation
\[%\begin{equation}\label{eq:2.10}
\varepsilon_{j-1}\varepsilon_j
b_{j-1}u_{j-1}-p_j(\lambda)u_{j}+b_{j}u_{j+1}=0\,\,\, j\in\dN,
\]%\end{equation}
obeying the initial conditions
\[%\begin{equation}\label{DiffEq}
 \begin{split}
    P_0(\lambda)&=1,\quad P_1(\lambda)=\frac{p_0(\lambda)}{b_0},\\
     Q_0(\lambda)&=0,\quad Q_1(\lambda)=\frac{\varepsilon_0}{b_0}.
\end{split}
\]%\end{equation}
According to~\eqref{SolMatr}, the $(j+1)$-th convergent of the
$P$-fraction~\eqref{ContF2} is equal to
\[
f_j:=T_{0}\circ T_{1}\circ\dots \circ
T_{j}(0)=-Q_{j+1}(\lambda)/P_{j+1}(\lambda).
\]
The relations~\eqref{SolMatr}, \eqref{W}, and \eqref{Wj} imply the
following relation~\cite{DD07}
\begin{equation}\label{Ostrogr2}
\varepsilon_jb_j(Q_{j+1}\left(\lambda)P_{j}(\lambda) -
        Q_{j}(\lambda) P_{j+1}(\lambda) \right)=1,\quad j\in\dZ_+.
\end{equation}

\subsection{Generalized Jacobi matrices}

 Let $p(\lambda)=p_{n}{\lambda}^{n}+\dots+p_{1}\lambda+p_{0}$
be a monic scalar real polynomial of degree $n$, i.e. $p_{n}=1$.
Let us associate with the polynomial $p$ its symmetrizator  $E_p$
and let the companion matrix $C_p$ be given by
\[%\begin{equation}\label{comp}
E_{p}=\begin{pmatrix}
p_{1}&\dots&p_{n}\\
\vdots&\sddots&\\
p_{n}&&{\bf 0}\\
\end{pmatrix},\quad
C_{p}=\begin{pmatrix}
0&\dots&0&-p_{0}\\
1&&{\bf 0}&-p_{1}\\
&\ddots&&\vdots\\
{\bf 0}&&1&-p_{n-1}\\
\end{pmatrix}.
\]%\end{equation}
As is known, $\det(\lambda-C_p)=p(\lambda)$ and the matrices $E_{p}$ and $C_{p}$ are related by 
\begin{equation}\label{cb}
C_{p}E_{p}=E_{p}C_{p}^{\top}.
\end{equation}

Now we are in a position to associate a tridiagonal block matrix with~\eqref{ContF2}.
\begin{definition}[\cite{DD}, \cite{KL79}]
Suppose we are given a $P$-fraction of the form~\eqref{ContF2}. 
Let $p_j$ be real scalar monic polynomials of degree ${k_j}$
\[
p_{j}(\lambda)={\lambda}^{k_{j}}+p_{k_{j}-1}^{(j)}{\lambda}^{k_{j}-1}+\dots+
p_{1}^{(j)}\lambda+p_{0}^{(j)},
\]
and let  $\e_{j}=\pm 1$,  $b_{j}>0$, $j\in\dZ_+$. The tridiagonal
block matrix
\begin{equation}\label{mJacobi}
H=\begin{pmatrix}
A_{0}   &\wB_{0}&       &{\bf 0}\\
B_{0}   &A_1    &\wB_{1}&\\
        &B_1    &A_{2} &\ddots\\
{\bf 0} &       &\ddots &\ddots\\
\end{pmatrix}
\end{equation}
where $A_{j}=C_{p_{j}}$ and $k_{j+1}\times k_{j}$ matrices $B_{j}$
and $k_{j}\times k_{j+1}$ matrices $\wB_{j}$ are given by
\begin{equation}\label{bblock}
B_{j}=\begin{pmatrix}
0&\dots&b_{j}\\
\hdotsfor{3}\\
0&\dots&0\\
\end{pmatrix},\,
\wB_{j}= \begin{pmatrix}
0&\dots&\wt b_{j}\\
\hdotsfor{3}\\
0&\dots&0\\
\end{pmatrix}, \quad
\wt b_{j}=\e_j\e_{j+1} b_j, \quad j\in\dZ_+,
\end{equation}
will be called a {\it generalized Jacobi matrix}  associated with
the $P$-fraction~\eqref{ContF2}.
\end{definition}
\begin{remark}
The papers~\cite{DD}, \cite{DD07}, and~\cite{KL79} are only
concerned with the case of generalized Jacobi matrices which are
finite rank perturbations of classical Jacobi matrices.
\end{remark}

Now, introducing the following shortened matrices
\begin{equation}\label{ShortMat}
{H}_{[0,j]}=
\begin{pmatrix}
A_{0}   &\wt{B}_{0}&       &\\
{B}_{0}   & {A}_{1}   & \ddots&\\
        & \ddots   & \ddots &\wt{B}_{j-1}\\
&       &{B}_{j-1} &{A}_j\\
\end{pmatrix},\,
{H}_{[1,j]}=
\begin{pmatrix}
{A}_{1}   &\wt{B}_{1}&       &\\
{B}_{1}   & {A}_{2}   & \ddots&\\
        & \ddots   & \ddots &{B}_{j-1}\\
&       &{B}_{j-1} &{A}_j\\
\end{pmatrix},
\end{equation}
where $j\in\dZ_+\cup\{\infty\}$, one can obtain the connection between the polynomials  ${P}_{j}$, ${Q}_{j}$ and the shortened Jacobi matrices
${H}_{[0,j]}$, $H_{[1,j]}$, respectively~\cite{DD}:
\[
 \begin{split}
P_{j}(\lambda)=(b_{0}\dots b_{j-1})^{-1}\det(\lambda-H_{[0,j-1]}),\\
Q_{j}(\lambda)=\varepsilon_0(b_{0}\dots b_{j-1})^{-1}\det(\lambda-H_{[1,j-1]}).
 \end{split}
\]

Furthermore, the following statement holds true.

\begin{proposition}[\cite{DD,DD07}]\label{prostota}
Let $j\in\dN$. Then
\begin{enumerate}
\item[i)] The polynomials $P_j$ and $P_{j+1}$  have no common zeros.
\item[ii)] The polynomials $Q_j$ and $Q_{j+1}$ have no common zeros.
\item[iii)] The polynomials $P_j$ and $Q_j$ have no common zeros.
\end{enumerate}
\end{proposition}

In what follows we are only interested in periodic generalized Jacobi matrices.

\begin{definition}\label{periodic_gjm}
Let $s\in\dN$. A generalized Jacobi matrix satisfying the
properties
\[
A_{js+k}=A_k,\quad B_{js+k}=B_{k},\quad\e_{js+k}=\e_k,\quad
j\in\dZ_+,\quad k\in\{0,\dots,s-1\}
\]
will be called an $s$-periodic generalized Jacobi matrix. 
\end{definition}

Let $\ell^2_{[0,\infty)}$ denote the Hilbert space of complex
square summable sequences $(w_0,w_1,\dots)$ equipped 
with the usual inner product. Evidently, any $s$-periodic generalized Jacobi matrix generates a
bounded linear operator in $\ell_{[0,+\infty)}^2$.
To see some more properties of periodic generalized Jacobi matrices
let us define the symmetric matrix $G$ by the equality
\begin{equation}\label{Gram}
G=\mbox{diag}(G_{0},G_1,\dots),\quad
G_{j}=\e_{j}E_{p_{j}}^{-1},\quad j\in\dZ_+.
\end{equation}
Clearly, the operator $G$ defined on $\ell^2_{[0,\infty)}$ is bounded and
self-adjoint. Moreover, $G^{-1}$ is a bounded linear operator in
$\ell^2_{[0,\infty)}$.

Let $\sH_{[0,\infty)}$ be a space of elements of $\ell^2_{[0,\infty)}$
 provided with the following indefinite inner product
\begin{equation}\label{metricInfty}
\left[x,y\right]=(Gx,y)_{\ell^2_{[0,\infty)}},\quad
x,y\in\ell^2_{[0,\infty)}.
\end{equation}

Let us recall~\cite{AI} that a pair $(\sH, [\cdot,\cdot])$
consisting of a Hilbert space $\sH$ and a sesquilinear form
$[\cdot,\cdot]$ on $\sH\times\sH$ is called {\it a space with
indefinite inner product}. A space with indefinite metric $(\sH,
[\cdot,\cdot])$ is called {\it a Krein space} if the indefinite
scalar product $[\cdot,\cdot]$ can be represented as follows
\[
[x,y]=(Jx,y)_{\sH}\quad x,y\in\sH,
\]
where the linear operator $J$ satisfies the following conditions
\[
J=J^{-1}=J^*.
\]
The operator $J$ is called {\it the fundamental symmetry}. So, one
can see that the space $\sH_{[0,\infty)}$ is the Krein space with
the fundamental symmetry $J=\sg G$ (see~\cite{AI} for details).
Moreover, the property~\eqref{cb} implies the following.

\begin{proposition}[\cite{De09,DD}]\label{SymGJM}
The $s$-periodic generalized Jacobi matrix defines a bounded
self-adjoint operator $H$ in the Krein space $\sH_{[0,\infty)}$, that is,
\[%\begin{equation}\label{hermitian}
\left[H x,y\right]=\left[x,H y\right]\quad x,y\in\sH_{[0,\infty)}.
\]%\end{equation}
\end{proposition}

The main tool for the spectral analysis of periodic generalized Jacobi
operators is the following matrix
\begin{equation}\label{monodromy}
T(\lambda):=\cW_{[0,s-1]}(\lambda)=\begin{pmatrix}
  -\e_{s-1}b_{s-1} Q_{s-1}(\lambda) & -Q_{s}(\lambda) \\
  \e_{s-1}b_{s-1} P_{s-1}(\lambda) & P_{s}(\lambda) \\
\end{pmatrix}.
\end{equation}
The matrix $T$ is called {\it the monodromy matrix}.
In fact, using the Floquet theory~\cite{De09} we can get the
description of the spectrum of the matrix $H$.
Indeed, let $w_1=w_1(\lambda)$ and $w_2=w_2(\lambda)$ be the roots of the
characteristic equation $\det(T(\lambda)-w)=0$. Introducing the
following notations
\[
\begin{split}
E&:=\{\lambda\in\dC: |w_1(\lambda)|=|w_2(\lambda)|\},\\
E_p&:=\{\lambda\in\dC:
P_{s-1}(\lambda)=0,\,|b_{s-1}Q_{s-1}(\lambda)|>|P_{s}(\lambda)|\}.
\end{split}
\]
we can give the description of spectra of periodic generalized Jacobi operators.
\begin{theorem}[\cite{De09}]\label{specper}
The spectrum $\sigma(H)$ of the $s$-periodic generalized Jacobi operator $H$ has
the form
\[
\sigma(H)=E\cup E_p,\quad \sigma_p(H)=E_p,
\]
where $\sigma_p(H)$ denotes the point spectrum of $H$, i.e. eigenvalues.
\end{theorem}

It should be mentioned that for the classical Jacobi matrices this result was obtained in~\cite{Ger}.

\begin{remark}[\cite{De09}]
One can also easily get another description of $E$
\[%\begin{equation}\label{h_ex_1}
E=\{\lambda\in\dC:
(P_s(\lambda)-\e_{s-1}b_sQ_{s-1}(\lambda))\in[-2,2]\}.
\]%\end{equation}
\end{remark}

\subsection{$m$-functions of periodic generalized Jacobi matrices}

Recall that the $m$-functions are one of the central tools in studying linear 
difference operators. 
We start with the definition of the $m$-function.
\begin{definition}\label{Weylfunction}
Let $H$ be a bounded generalized Jacobi matrix. The function $m$ defined by
\begin{equation}\label{Wfexf}
m(\lambda)=[(H-\lambda)^{-1}e,e],\quad e=(1,0,0,\dots)^{\top}
\end{equation}
 is called the $m$-function (or, the Weyl function) of the operator $H$.
\end{definition}
Next, by using the Frobenius formula, one can see that the $m$-function
$m$ and the $m$-function
\[
m_{[1,\infty)}(\lambda)=[(H_{[1,\infty)}-\lambda)^{-1}e,e]=
(G_{[1,\infty)}(H_{[1,\infty)}-\lambda)^{-1}e,e)
\]
of $H_{[1,\infty)}$ (for the definition of $H_{[1,\infty)}$ and $G_{[1,\infty)}$ see~\eqref{ShortMat}) are related by the equality~\cite{DD}
\begin{equation}\label{Riccatiinfty}
m(\lambda)=\frac{-\varepsilon_0}
    {p_0(\lambda)+\varepsilon_0b_0^2m_{[1,\infty)}(\lambda)},\quad |\lambda|>\Vert H\Vert\ge\Vert H_{[1,\infty)}\Vert.
\end{equation}

 The standard technique of dealing with periodic continued fractions enables us to calculate the $m$-function explicitly.
\begin{proposition}\label{mfPer}
 Let $H$ be an $s$-periodic generalized Jacobi matrix. Then its $m$-function has the following form
\begin{equation}\label{forMFper}
m(\lambda)=\frac{-(P_s(\lambda)+\varepsilon_{s-1}b_{s-1}Q_{s-1}(\lambda))+
\sqrt{(P_s(\lambda)-\varepsilon_{s-1}b_{s-1}Q_{s-1}(\lambda))^2-4}}{2\varepsilon_{s-1}b_{s-1}P_{s-1}(\lambda)},
\end{equation}
where the cut is $E$ and the branch is chosen in such a way that $m(\lambda)\to 0$ as $\lambda\to\infty$.
\end{proposition}
\begin{proof}
Due to~\eqref{Riccatiinfty} and the $s$-periodicity of $H$, we have that
\[
m(\lambda)=
-\frac{\varepsilon_0}{p_0(\lambda)}
\begin{array}{l} \\ - \end{array}
\frac{\varepsilon_0\varepsilon_1b_0^2}{p_1(\lambda)}
\begin{array}{ccc} \\ - & \cdots & -\end{array}
\frac{\varepsilon_{s-1}\varepsilon_{s-2}b_{s-2}^2}
{p_{s-1}(\lambda)+\varepsilon_{s-1}b_{s-1}^2m(\lambda)}.
\]
According to~\eqref{W} and~\eqref{monodromy} the latter relation takes the form 
\begin{equation}\label{mfH1}
m(\lambda)=-\frac{\varepsilon_{s-1}b_{s-1}Q_{s-1}(\lambda)m(\lambda)+Q_{s}(\lambda)}
{\varepsilon_{s-1}b_{s-1}P_{s-1}(\lambda)m(\lambda)+P_{s}(\lambda)}.
\end{equation}
Next,~\eqref{mfH1} can be rewritten as follows
\begin{equation}\label{mfH2}
\varepsilon_{s-1}b_{s-1}P_{s-1}(\lambda)m^2(\lambda)
+(P_s(\lambda)+\varepsilon_{s-1}b_{s-1}Q_{s-1}(\lambda))m(\lambda)+Q_{s}(\lambda)=0.
\end{equation}
Now,~\eqref{forMFper} follows from~\eqref{mfH2} and~\eqref{Ostrogr2}.
The choice of the cut and the branch is due to Definition~\ref{Weylfunction} and 
Theorem~\ref{specper}.
\end{proof}
\begin{theorem}[cf.~\cite{De09}]
 The $m$-function of the periodic generalized Jacobi matrix $H$ determines $H$ uniquely.
\end{theorem}
\begin{proof}
 In fact,
consecutive applications of the relation~\eqref{Riccatiinfty}
lead to the continued fraction~\eqref{ContF2}, which determines $H$ uniquely. 
\end{proof}

\section{Inverse problems}

\subsection{The monodromy matrices} Note that knowing the monodromy matrix $T$ gives 
the complete information about the spectrum of the corresponding $s$-periodic generalized Jacobi matrix $H$.
Moreover, it is clear that one can reconstruct the $s$-periodic generalized Jacobi matrix $H$ by its monodromy matrix
$T=\cW_{[0,s-1]}$ since by the construction we have that
\[
 -\frac{Q_{s}(\lambda)}{P_s(\lambda)}=
-\frac{\varepsilon_0}{p_0(\lambda)}
\begin{array}{l} \\ - \end{array}
\frac{\varepsilon_0\varepsilon_1b_0^2}{p_1(\lambda)}
\begin{array}{ccc} \\ - & \cdots & -\end{array}
\frac{\varepsilon_{s-2}\varepsilon_{s-1}b_{s-2}^2}{p_{s-1}(\lambda)}
\]
and $b_{s-1}$ can be extracted  from $T$ directly.
Now, one of the most natural questions is the following: which matrices can be the monodromy matrices?
To answer this question let us note first that the monodromy matrix 
\[
 T(\lambda)=
\begin{pmatrix}
t_{11}(\lambda)& t_{12}(\lambda)\\
t_{21}(\lambda)& t_{22}(\lambda)\\
\end{pmatrix}
\]
satisfies the following properties
\begin{enumerate}
 \item[(T1)] $t_{11}$, $t_{12}$, $t_{21}$, $t_{22}$ are polynomials with real coefficients such that
\[
 t_{11}(\lambda)t_{22}(\lambda)-t_{12}(\lambda)t_{21}(\lambda)=1, \quad \lambda\in\dR;
\]
\item[(T2)] $\deg t_{12}<\deg t_{22}$, $\deg t_{21}<\deg t_{22}$;
\item[(T3)] the absolute values of the leading coefficients of $t_{21}$ and $t_{22}$ are equal and
the leading coefficient of $t_{22}$ is positive.
\end{enumerate}
These properties give us a tip to introduce the following definition.
\begin{definition}\label{admis}
 The $2\times 2$-matrix polynomial $T$ is called admissible if it satisfies the conditions (T1)-(T3).
\end{definition}
\begin{remark}\label{rem32}
 Consider the following $2\times 2$-matrix
\[
 J=\begin{pmatrix}
0 & -i\\
i & 0\\
   \end{pmatrix}.
\]
It is said that $2\times 2$ matrix polynomial is $J$-unitary on the real line $\dR$
if the following equality holds true
\[
 T(\lambda)JT^*(\overline{\lambda})=J, \quad \lambda\in\dR.
\]
Clearly, every admissible matrix polynomial is $J$-unitary. Moreover, every
$J$-unitary $2\times 2$-matrix polynomial on the real line can be normalized 
by multiplications from left and right by constant $J$-unitary matrices to an admissible matrix
(for more details see~\cite[Section~3.7]{DD07}).
\end{remark}

The answer to the above-stated question is given by the following statement. 
\begin{theorem}[cf.~\cite{DD07}]\label{InvMon}
Every admissible matrix polynomial $T$ is the monodromy matrix of a unique periodic generalized Jacobi matrix.
 \end{theorem}
\begin{proof}
Actually, this statement can be proved by re-examining~\cite[Section~3.7]{DD07}.
First, note that the admissible matrix $T$ generates the shortened generalized Jacobi matrix $H_{[0,s-1]}$.
Indeed, let us consider the $P$-fraction expansion 
\[
 \frac{t_{12}(\lambda)}{t_{22}(\lambda)}=
-\frac{\varepsilon_0}{p_0(\lambda)}
\begin{array}{l} \\ - \end{array}
\frac{\varepsilon_0\varepsilon_1b_0^2}{p_1(\lambda)}
\begin{array}{ccc} \\ - & \cdots & -\end{array}
\frac{\varepsilon_{s-2}\varepsilon_{s-1}b_{s-2}^2}{p_{s-1}(\lambda)},
\]
which, according to the definition, generates the shortened generalized Jacobi matrix $H_{[0,s-1]}$.
So, one can recover the sequences $P_0$, $P_1$, \dots, $P_{s-1}$, $b_{s-1}P_s$ and 
$Q_0$, $Q_1$, \dots, $Q_{s-1}$, $b_{s-1}Q_s$ such that $t_{12}/t_{22}=-Q_{s}/P_{s}$.
Choosing $b_{s-1}>0$  to have the equalities 
$t_{12}=-Q_s$ and $t_{22}=P_s$ we conclude that we know the matrix
\[
\cW_{[0,s-1]}(\lambda)=\begin{pmatrix}
  -\e_{s-1}b_{s-1} Q_{s-1}(\lambda) & -Q_{s}(\lambda) \\
  \e_{s-1}b_{s-1} P_{s-1}(\lambda) & P_{s}(\lambda) \\
\end{pmatrix}
\]
and the number $b_{s-1}$.
Also, the periodicity implies that $\varepsilon_{0}=\varepsilon_{s}$ and, so, one can 
determine $\wb_{s-1}=\varepsilon_{s-1}\varepsilon_{s}b_{s-1}$. Thus, the admissible matrix $T$
gives rise to the following $s$-periodic generalized Jacobi matrix
\[
H=
\begin{pmatrix}
H_{[0,s-1]}   &\wB_{s-1}&       &{\bf 0}\\
B_{s-1}   & H_{[0,s-1]}    &\wB_{s-1}&\\
        &B_{s-1}    &H_{[0,s-1]} &\ddots\\
{\bf 0} &       &\ddots &\ddots\\
\end{pmatrix}.
\]
Finally, it follows from (T1)-(T3), Proposition~\ref{prostota}, and~\eqref{Ostrogr2} that 
the monodromy matrix $\cW_{[0,s-1]}$ of $H$ coincides with $T$
(for more details see the proof of~\cite[Theorem~3.23]{DD07}).
\end{proof}
\begin{remark}
 Due to Remark~\ref{rem32} and formula~\eqref{W}, Theorem~\ref{InvMon} can be applied
for getting factorization results for $J$-unitary matrix polynomials~\cite{ADL04},~\cite{DD07}.
\end{remark}

\subsection{The Abel criterion}  
Formula~\eqref{forMFper} leads to the conclusion that the $m$-function of the periodic
generalized Jacobi matrix has the following algebraic form
\begin{equation}\label{hyperell}
\varphi(\lambda)=\frac{\sqrt{R(\lambda)}-U(\lambda)}{V(\lambda)},
\end{equation}
where $R$, $U$, and $V$ are polynomials such that
\begin{equation}\label{degP}
 \deg R=2n,\quad \deg U=n, \quad \deg V<\deg U.
\end{equation}
 Besides, the $m$-function  has the property that $\varphi(\lambda)\to 0$ as $\lambda\to\infty$.

Clearly, the function $\varphi$ can be also represented by means of some other choice of polynomials
$R$, $U$, and $V$. So, without loss of generality we always assume that the polynomials
$R$, $U$, and $V$ are chosen in such a way that the polynomial $V$ has the minimal degree between all
such polynomials in the representation~\eqref{hyperell} of the given function $\varphi$.
Now the natural question is to ask under what conditions on $R$, $U$, and $V$ a function of
the form~\eqref{hyperell} can be the $m$-function of a periodic generalized Jacobi matrix?
One can find the answer in the following statement.
\begin{theorem}\label{PellAbelTh}
 Let $R$, $U$, and $V$ be polynomials subject to~\eqref{degP}. Let us also suppose that one of the  branches of the function $\varphi=(\sqrt{R}-U)/V$ possesses the property
$\varphi(\lambda)\to 0$ as $\lambda\to\infty$. Then there is a cut such that the corresponding branch of
$\varphi$ is the $m$-function of a periodic generalized Jacobi matrix if and only if there exist
real polynomials $X$, $Y$, $Z$  satisfying the following relations
\begin{eqnarray}
X^2-RY^2=1,\label{PellAbel1}\\
(U^2-R)Y=VZ.\label{PellAbel2}
\end{eqnarray}
\end{theorem}
\begin{proof}{\it The "if" part.} Let $\varphi=(\sqrt{R}-U)/V$ be the $m$-function
of a periodic generalized Jacobi matrix. Then $\varphi$ satisfies the following equation
\begin{equation}\label{draft63}
\varepsilon_{s-1}b_{s-1}P_{s-1}(\lambda)\varphi^2(\lambda)
+(P_s(\lambda)+\varepsilon_{s-1}b_{s-1}Q_{s-1}(\lambda))\varphi(\lambda)+Q_{s}(\lambda)=0.
\end{equation}
Substituting $\varphi=(\sqrt{R}-U)/V$ to~\eqref{draft63} one gets
\begin{equation}\label{draft64}
 \begin{split}
 \sqrt{R}\Big(-2\varepsilon_{s-1}b_{s-1}P_{s-1}U+
 P_sV+\varepsilon_{s-1}b_{s-1}Q_{s-1}V\Big)+\\
 +\varepsilon_{s-1}b_{s-1}P_{s-1}(R^2+U^2)-
 (P_s+\varepsilon_{s-1}b_{s-1}Q_{s-1})UV+
 Q_sV^2=0.
 \end{split}
\end{equation}
Observe that the equality $\alpha\sqrt{R}+\beta=0$, where $\alpha$ and $\beta$ are some polynomials, implies 
$\alpha=0$ and $\beta=0$. So, it follows from~\eqref{draft64} that
\begin{eqnarray}
P_sV&=&2\varepsilon_{s-1}b_{s-1}P_{s-1}U-\varepsilon_{s-1}b_{s-1}Q_{s-1}V,\label{draft65}\\
  Q_sV^2&=&(P_s+\varepsilon_{s-1}b_{s-1}Q_{s-1})UV-\varepsilon_{s-1}b_{s-1}P_{s-1}(R^2+U^2). \label{draft66}
\end{eqnarray}
Combining~\eqref{draft65} and~\eqref{draft66} gives
\begin{equation}\label{draft67}
Q_sV^2=\varepsilon_{s-1}b_{s-1}P_{s-1}(U^2-R).
\end{equation}
It is easy to see from~\eqref{forMFper} as well as from~\eqref{mfH2} that $P_{s-1}=\widetilde{P}_{s-1}V$.
Taking this observation into account one can rewrite~\eqref{draft65} and~\eqref{draft67} as follows
\begin{eqnarray}
P_s&=&2\varepsilon_{s-1}b_{s-1}\widetilde{P}_{s-1}U-\varepsilon_{s-1}b_{s-1}Q_{s-1},\label{draft68}\\
  Q_sV&=&\varepsilon_{s-1}b_{s-1}\widetilde{P}_{s-1}(U^2-R). \label{draft69}
\end{eqnarray}
Now, by setting $Z=Q_{s}$ and $Y=\varepsilon_{s-1}b_{s-1}\widetilde{P}_{s-1}$ the relation~\eqref{draft69}
becomes~\eqref{PellAbel2}.

In order to get~\eqref{PellAbel1}, let us multiply formula~\eqref{Ostrogr2} by $V$
\[
\varepsilon_{s-1}b_{s-1}\left(VQ_{s}P_{s-1}-VQ_{s-1} P_{s} \right)=1.
\]
 Substituting~\eqref{draft68} and~\eqref{draft69} to the latter equality we arrive at 
\[
\varepsilon_{s-1}b_{s-1}\Big(\varepsilon_{s-1}b_{s-1}\widetilde{P}_{s-1}^2V(U^2-R)
-\varepsilon_{s-1}b_{s-1}V(2U\widetilde{P}_{s-1}-Q_{s-1})Q_{s-1}\Big)=V.
\]
Further calculations lead to the relation
\[
 b_{s-1}^2\Big(\widetilde{P}_{s-1}^2(U^2-R)-2U\widetilde{P}_{s-1}Q_{s-1}+Q_{s-1}^2\Big)=1,
\]
which can be rewritten as follows
\begin{equation}\label{draft611}
 b_{s-1}^2(\widetilde{P}_{s-1}U-Q_{s-1})^2-b_{s-1}^2R\widetilde{P}_{s-1}^2=1.
\end{equation}
Setting $X=\varepsilon_{s-1}b_{s-1}(\widetilde{P}_{s-1}U-Q_{s-1})$ and recalling that
 $Y=\varepsilon_{s-1}b_{s-1}\widetilde{P}_{s-1}$ 
one sees that~\eqref{draft611} is exactly~\eqref{PellAbel1}.

{\it The "only if" part.} Suppose that the system~\eqref{PellAbel1},~\eqref{PellAbel2} is satisfied for 
some polynomials $X$, $Y$, and $Z$. Then the ''if'' part gives us a hint to introduce the following 
polynomials
\[
t_{11}=X-YU,\quad t_{12}=YV,\quad t_{21}=-Z,\quad t_{22}=X+YU,
\]
where the signs of the leading coefficients of $X$, $Y$, and $Z$ are taken to have
the leading coefficient $t_{22}$ positive and $\deg t_{11}<n$.
Next it is easy to check 
\[
\begin{split}
\det T&=t_{11}t_{22}-t_{12}t_{21}=X^2-Y^2U^2+ZYV=\\
&=X^2-(RY^2-ZYV)+ZYV=X^2-RY^2=1.
\end{split}
\]
It remains to note that 
\[
\deg t_{21}=\deg Z=\deg (U^2-R)Y/V<\deg X=\deg t_{22}
\]
in order to see that the matrix $T=(t_{ij})_{i,j=1}^2$ is admissible.
The rest follows from Theorem~\ref{InvMon}.
\end{proof}
\begin{remark}
 If $V=1$ then Theorem~\ref{PellAbelTh} reduces to the classical Abel criterion of the periodicity of continued fractions
representing $\sqrt{R}$.
\end{remark}

\begin{remark}
 Every solvable system~\eqref{PellAbel1},~\eqref{PellAbel2} leads to an admissible matrix polynomial. The converse is also true.
In turn, every admissible matrix polynomial can be constructed by~\eqref{Wj},~\eqref{W}.
\end{remark}

\noindent{\bf Acknowledgments}. This work was mainly done several years ago when I was
 a PhD student at the Donetsk National University. I would like to express  my gratitude to
my scientific adviser Professor V.A.~Derkach for the numerous discussions and helpful comments.
I am also deeply indebted to Professors V.P.~Burskii and A.S.~Zhedanov who involved me in studying
the papers~\cite{Mal01} and~\cite{SYu}.

\end{document}